\documentclass[11pt,reqno]{amsart}
\usepackage{amscd,amssymb,amsmath,amsthm,amsfonts}
\usepackage{graphicx}

\oddsidemargin=0.3in \evensidemargin=0.3in \theoremstyle{plain}
\newtheorem{thm}[subsection]{Theorem}
\newtheorem{lemma}[subsection]{Lemma}
\newtheorem{proposition}[subsection]{Proposition}
\newtheorem{cor}[subsection]{Corollary}

\newtheorem{rk}[subsection]{Remark}
\newtheorem{defn}[subsection]{Definition}
\newtheorem{ex}[subsection]{Example}

\numberwithin{equation}{section} 

\newcommand{\R}{\mathbb{R}}

\begin{document}
\title [On nilpotent index and dibaricity of evolution algebras]{On nilpotent index and dibaricity of evolution algebras}

\author {J.M. Casas, M. Ladra, B.A. Omirov, U.A. Rozikov}

\address{J.M.\ Casas\\ Department of Applied Mathematics I, E.E. Forestal, Pontevedra, University of Vigo, 36005, Spain.}
 \email {jmcasas@uvigo.es}
 \address{M.\ Ladra\\ Department of Algebra, University of Santiago de Compostela, 15782, Spain.}
 \email {manuel.ladra@usc.es}
 \address{B.A.\ Omirov and U.A.\ Rozikov\\ Institute of mathematics,
Tashkent, 100125, Uzbekistan.} \email {omirovb@mail.ru,
rozikovu@yandex.ru}



\begin{abstract} An evolution algebra corresponds to a quadratic matrix $A$ of structural constants. It is known
the equivalence between nil, right nilpotent evolution
algebras and evolution algebras which are defined by upper
triangular matrices $A$. We establish a criterion for an $n$-dimensional nilpotent evolution algebra
to be with maximal nilpotent index $2^{n-1}+1$. We give the classification of finite-dimensional complex evolution algebras with maximal nilpotent index. Moreover, for any $s=1,\dots,n-1$ we construct a wide class of $n$-dimensional evolution algebras with
nilpotent index $2^{n-s}+1$. We show that nilpotent evolution algebras are not dibaric  and establish a criterion for two-dimensional real evolution algebras to be dibaric. 
\end{abstract}

\maketitle

\section{Introduction}

{\bf A motivation for replacement of the first version:} {\sl The first version of our paper was in arXiv:1208.1709, which was also 
published in \cite{LAA}. The paper \cite{LAA} was reviewed (for MathSciNet) by J. Carlos Guti\'errez Fern\'andez. The reviewer mentioned that ``{\it...there are a few mistakes in the article. Lemma 3.2 asserts that $E^i, i\geq 1$, is an evolution subalgebra of $E$, in the sense that $E^i$ admits a natural basis which can be extended to a natural basis of $E$. The reviewer thinks this statement is false. This claim is not true because the example given in Remark 3.14 satisfies $E^3\ne 0$ and $E^4=0$. Also, the example of Remark 3.14 for $f=1$ contradicts the claim in Proposition 3.11}". Here we note that in Lemma 3.2 we are not asserting but {\bf assuming} that $E^i, i\geq 1$, is an evolution subalgebra of $E$.
We have slightly corrected the conditions of the Proposition 3.11. So mistakes mentioned by the reviewer are corrected in this version.}\\

In the book \cite{Tian1} the foundations of evolution algebras
are developed. Evolution algebras have many
connections with other mathematical fields including graph theory, group theory,
Markov chains, dynamic systems, knot theory, 3-manifolds and the study
of the Riemann-zeta functions \cite{Omirov,CLOR,LOR,Tian2,Tian1}.

Nilpotent algebra is an algebra for which there is a natural number $k$ such that any product of $k$ elements of the algebra is zero. If there is a non-zero product of $k-1$ elements, then $k$ is called the index of nilpotency  of the algebra. Examples of nilpotent algebras are: an algebra with zero multiplication; direct sums of nilpotent algebras, the nilpotent indices of which are uniformly bounded; and the tensor product of two algebras, one of which is nilpotent.
Nilpotent subalgebras that coincide with their normalizer (Cartan subalgebras) play an essential role in the classification of simple Lie algebras of finite dimension.

The algebraic notions like nilpotency, right nilpotency and solvability might be interpreted in a biological way as a various types of vanishing (``deaths'') populations.

The structural constants of an evolution
algebra are given by a quadratic matrix $A$ (see Section \ref{S:prel}). In \cite{CLOR} the equivalence between nil, right nilpotent evolution
algebras and evolution algebras which are defined by upper
triangular matrices $A$ is proved, and the classification of 2-dimensional complex
evolution algebras is obtained.

In \cite{CG} the derivations of $n$-dimensional complex evolution algebras, depending on the rank of the matrix $A$, are studied. For an evolution algebra with non-singular matriz it is proved that the space of derivations is zero. The spaces of derivations for evolution algebras with matrices of rank $n-1$ are described.

The paper \cite{Omirov} is devoted to the study of finite-dimensional complex evolution algebras. The class of evolution algebras isomorphic to evolution algebras with
Jordan form matrices of structural constants is described. For finite-dimensional complex evolution algebras
the criteria of nilpotency is established in terms of the properties of the corresponding
matrices. Moreover, it is proved that for nilpotent $n$-dimensional complex evolution
algebras the possible maximal nilpotency index is $2^{n-1}+1$. The criteria of planarity
for finite graphs is formulated by means of evolution algebras defined by graphs.

In \cite{ly} an evolution algebra $E$ associated to the
free population is introduced and using this non-associative baric
algebra, many results are obtained in an explicit form, e.g. the
explicit description of stationary quadratic operators and the
explicit solutions of a non-linear evolutionary equation in the
absence of selection, as well as general theorems on convergence
to equilibrium in the presence of a selection.

Dibaric algebras have not non-zero homomorphisms to the set of the real numbers.
In \cite{LOR} a concept of bq-homomorphism (which is given by two linear maps $f, g$
of the algebra to the set of the real numbers) is introduced and it is shown that an
 algebra is dibaric if and only if it admits a non-zero bq-homomorphism.

In the study of any class of algebras, it is important to describe, up
to isomorphism, at least algebras of lower dimensions because such
description gives examples to establish or reject certain
conjectures. In this way in \cite{M} and \cite{Umlauf}, the
classifications of associative and nilpotent Lie algebras of low
dimensions were given.

In this paper we continue the study of algebraic properties of evolution algebras. The
paper is organized as follows. In Section \ref{S:prel} we give some preliminaries. In Section \ref{S:nil}
we establish a criterion for an $n$-dimensional nilpotent evolution algebra
to be of maximal nilpotent index $2^{n-1}+1$. Since these algebras have maximal index of right
 nilpotency and maximal index of solvability too, then we might say that among vanishing populations,
  these are the latest vanishing populations.

We give the classification of finite-dimensional complex evolution algebras with maximal nilpotent index. Moreover, for any $s=1,\dots,n-1$ we construct a wide class of $n$-dimensional evolution algebras with
nilpotent index $2^{n-s}+1$. Section \ref{S:dib} is devoted to the dibaricity of evolution algebras. We show that nilpotent real evolution algebras are not dibaric and establish a criterion for two-dimensional real evolution algebras to be dibaric.
\section{Preliminaries} \label{S:prel}

\subsection*{Evolution algebras} Let $(E,\cdot)$ be an algebra over a field $K$. If it admits a basis
$\{e_1,e_2,\dots\}$, such that
\[
e_i \cdot e_j=
\begin{cases}
0, &\text{if \ $i\ne j$;}\\
\displaystyle \sum_{k}a_{ik}e_k, &\text{if \ $i=j$,}
\end{cases}
\]
then this algebra is called an {\it evolution algebra} \cite{Tian1}. The basis is called a natural basis.  We denote by
$A=(a_{ij})$ the matrix of the structural constants of the evolution
algebra $E$.

The following properties are known \cite{Tian1}:
\begin{enumerate}
  \item Evolution algebras are not associative, in general.
  \item Evolution algebras are commutative, flexible.
  \item Evolution algebras are not power-associative, in general.
  \item The direct sum of evolution algebras is also an evolution algebra.
  \item The Kronecker product of evolution algebras is an evolution algebra.
\end{enumerate}

\begin{defn} Let $E$ be an evolution algebra, and $E_1$ be a subspace of $E$.
If $E_1$ has a natural basis $\{e_i \mid i\in \Lambda_1\}$, which can be extended to a natural
basis $\{e_j \mid j\in \Lambda\}$ of $E$, $E_1$ is called an evolution subalgebra, where $\Lambda_1$ and $\Lambda$
are index sets and $\Lambda_1$ is a subset of $\Lambda$.
\end{defn}

\begin{defn} An element $a$ of an evolution algebra $E$  is called nil if
there exists $n(a) \in \mathbb{N}$ such that $\big(\cdots
\underbrace{\big((a\cdot a)\cdot a\big)\cdots a}_{n(a)}\big)=0$.
An evolution algebra $E$ is called nil if every  element of the algebra
is nil.
\end{defn}

For an evolution algebra $E$ we introduce the following sequences, \  $k \geq$ 1,
\[
E^{[1]}=E^{<1>}=E, \qquad  E^{[k+1]}=E^{[k]}E^{[k]}, \qquad  E^{<k+1>}=E^{<k>}E,
\]
\begin{equation}\label{EE}
E^k=\sum_{i=1}^{k-1}E^iE^{k-i}.
\end{equation}
The following inclusions hold
\[E^{<k>}\subseteq E^k, \qquad  E^{[k+1]}\subseteq E^{2^k}.\]
Since $E$ is a commutative algebra we obtain
\[E^k=\sum_{i=1}^{\lfloor k/2\rfloor}E^iE^{k-i},\]
where $\lfloor x\rfloor$ denotes the integer part of $x$.
\begin{defn} An evolution algebra $E$ is called right nilpotent if there
exists some $s\in \mathbb{N}$ such that $E^{<s>}=0$. The smallest $s$ such that
 $E^{<s>}=0$ is called the index of right nilpotency.
\end{defn}

\begin{defn} An evolution algebra $E$ is called  nilpotent if there
exists some $n\in \mathbb{N}$ such that $E^n=0$. The smallest $n$ such that
 $E^n=0$ is called the index of  nilpotency.
\end{defn}

In \cite{Omirov}, it is proved that the notions of nilpotent and right
nilpotent are equivalent.

\begin{defn} An algebra $\mathbf A$ is called {\it solvable} if there exits some $t \in \mathbb{N}$ such that $\mathbf A^{[t]}=0$.
The smallest $t$ such that $\mathbf A^{[t]}=0$ is called the index of solvability.
\end{defn}

\subsection*{Dibaric algebras} A \emph{character} for an algebra $\mathbf A$ is a non-zero multiplicative
linear form on $\mathbf A$, that is, a non-zero algebra homomorphism from $\mathbf A$
to $\R$ \cite{ly}. A pair $(\mathbf A, \sigma)$ consisting of an algebra $\mathbf A$ and a character
$\sigma$ on $\mathbf A$ is called a \emph{baric algebra}.

As usual,  the algebras considered in mathematical biology are not baric.

\begin{defn}[\cite{m,w}] Let $\mathfrak A=\langle w, m \rangle_\R$ denote
a two-dimensional commutative algebra over $\R$ with multiplicative
table
\[w^2=m^2=0, \qquad wm= \frac{1}{2} (w+m)\, .\]
Then  $\mathfrak A$ is called the sex differentiation algebra.
\end{defn}

 It is clear that
$\mathfrak A^2=\langle w+m \rangle_\R$ is an ideal of $\mathfrak A$
which is isomorphic to the field $\R$. Hence the algebra $\mathfrak
A^2$ is a baric algebra.

\begin{defn}[\cite{m}] An algebra $\mathbf A$ is called dibaric if it admits a
homomorphism onto the sex differentiation algebra $\mathfrak A$.
 \end{defn}

\section{Nilpotent evolution algebras} \label{S:nil}

In \cite{CLOR} it is proved that the notions of nil and right nilpotency
are equivalent for evolution algebras. Moreover, the matrix $A$ of
structural constants for such algebras has upper (or lower, up to
permutation of basis elements of the algebra) triangular form.

Let evolution algebra $E$ be a right nilpotent algebra, then it is
evident that $E$ is a nil algebra. Therefore for the related
matrix $A=\left(a_{ij}\right)_{i,j=1}^n$, we have
\[a_{i_1i_2}a_{i_2i_3}\dots a_{i_ki_1}=0,\]
for any $k\in \{1, 2, \dots,
n\}$ and arbitrary $i_1, i_2, \dots, i_k\in \{1,2, \dots, n\}$ with
$i_p\ne i_q$ for $p\ne q$ \cite{CLOR}.

The following results are known:

\begin{thm}[\cite{CLOR}]   The following statements are equivalent for
an $n$-dimensional evolution algebra $E$:

\begin{itemize}
\item[(a)] The matrix corresponding to $E$ can be written as
\[
\widehat{A}=
  \begin{pmatrix}
 0 & a_{12} & a_{13} &\dots &a_{1n} \\[1.5mm]
 0 & 0 & a_{23} &\dots &a_{2n} \\[1.5mm]
0 & 0 & 0 &\dots &a_{3n} \\[1.5mm]
\vdots & \vdots & \vdots &\cdots & \vdots \\[1.5mm]
0 & 0 & 0 &\cdots &0
\end{pmatrix};
\]

\item[(b)] $E$ is a right nilpotent algebra;

\item[(c)] $E$ is a nil algebra.
\end{itemize}
\end{thm}

\begin{lemma}\label{ln} Let $E$ be a finite-dimensional evolution algebra and $E^j$, $j\geq 1$,
the evolution subalgebras of $E$ defined in \eqref{EE}.
Then
\begin{equation}\label{k}
E^{2^k+i}=E^{2^{k+1}}, \quad  i=1,\dots,2^k, \ \ k=0,1,\dots.
\end{equation}
\end{lemma}
\proof We shall use mathematical induction. We have
$E^1=E$, $E^2=EE$, and for $k=1$,
\[ E^3=EE^2=E^2E^2, \qquad \ E^4=EE^3+E^2E^2=E^2E^2=E^3.\]
Assume for $k$ the equalities \eqref{k} are true. We shall prove for $k+1$.
Using $E^{i+j}\subset E^i$, $E^{i+j}E^i=E^{i+j}E^{i+j}$ and assumptions of the induction we get
\begin{align*}
E^{2^{k+1}+i}&=\sum_{j=1}^{2^k+\lfloor i/2\rfloor}E^jE^{2^{k+1}+i-j}=\sum_{j=1}^{2^k+\lfloor i/2\rfloor}E^jE^{2^k+i}\\
{}&=\sum_{j=1}^{2^k+\lfloor i/2\rfloor}E^jE^{2^{k+1}}=EE^{2^{k+1}}=E^{2^{k+1}}E^{2^{k+1}}. \qedhere
\end{align*}
\endproof
\begin{lemma}\label{lb} If for an evolution algebra $E$ there exists $s\in \mathbb{N}$ such that $E^{2^s+1}=E^{2^{s+1}+1}$,
then $E^k=E^{2^s+1}$ for any $k=2^s+1, 2^s+2, \dots, 2^{s+2}+1$.
\end{lemma}
\proof We have
\[E^{2^s+1}\supseteq E^{2^{s}+2}\supseteq \dots \supseteq E^{2^{s+1}+1}.\]
Hence by condition of the lemma we get
 $E^k=E^{2^s+1}$ for any $k=2^s+1, 2^s+2, \dots, 2^{s+1}+1$.
 It remains to prove the equality for $k=2^{s+1}+i+1$, $i=1,\dots, 2^{s+1}$.
 We have
\[E^{2^{s+1}+1}=\sum_{j=1}^{2^s}E^jE^{2^{s+1}-j+1}=EE^{2^s+1}=E^{2^s+1}.\]
 For $i=1$ using the above obtained equalities we get
\[E^{2^{s+1}+2}=\sum_{j=1}^{2^s+1}E^jE^{2^{s+1}-j+2}=EE^{2^s+1}=E^{2^s+1}.\]
 Now assume the assertion is true for $i$ and we shall show it for $i+1$.
\[E^{2^{s+1}+i+2}=\sum_{j=1}^{2^s+1+\lfloor i/2\rfloor}E^jE^{2^{s+1}+i-j+2}.\]
Since $2^{s}+1<2^{s+1}+i-j+2\leq 2^{s+1}+i+1$ for any $j=1,2,\dots 2^s+1+\lfloor i/2\rfloor$, using the assumption of the induction, we get
\[E^{2^{s+1}+i+2}=\sum_{j=1}^{2^s+1+\lfloor i/2\rfloor}E^jE^{2^{s}+1}=EE^{2^s+1}=E^{2^s+1}. \qedhere \]
\endproof
From this lemma we get the following
\begin{cor} If for an evolution algebra $E$ there exists $s\in \mathbb{N}$ such that $E^{2^s+1}=E^{2^{s+1}+1}$,
then $E^k=E^{2^s+1}$ for any $k\geq 2^s+1$.
\end{cor}
\proof If the condition of Lemma \ref{lb} is satisfied for $s$, then it is satisfied for $s+1$. So, iterating
the lemma we get $E^k=E^{2^s+1}$ for any $k\geq 2^s+1$.
\endproof
From this corollary it follows that an evolution algebra $E$ satisfying the condition of Lemma \ref{ln} is not nilpotent.

The following is an example satisfying the condition of Lemma \ref{lb}.
\begin{ex} Fix some $r\in \{2,3,\dots, n-1\}$ and consider the evolution algebra with the multiplication table
\[e^2_i=e_{i+1}, \ \  i=1,\dots,r-1; \qquad  e^2_i=e_r,  \ \ i=r,\dots,n.\] It is easy to see that this algebra satisfies the condition of Lemma \ref{lb} for some $s\geq r$. In this case, $E^k=\{e_r\}$ for all $k\geq 2^s+1$.
\end{ex}

\begin{thm} An $n$-dimensional nilpotent evolution algebra $E$ has maximal nilpotent index, $2^{n-1}+1$, if and only if
\[a_{12}a_{23}\dots a_{n-1,n}\ne 0.\]
\end{thm}
\proof {\sl Necessity.} Assume $a_{12}a_{23}\dots a_{n-1,n}=0$ then $\dim E^2\leq n-2$. Since $E$ is nilpotent, by Lemma \ref{ln}, for any $k$
we have $E^{2^k+1}\supsetneqq E^{2^{k+1}+1}$. Consequently,  $\dim E^{2^k+1}\leq n-2-k$.
Hence $E^{2^{n-2}+1}=0$, i.e., $E$ has not maximal nilpotent index.

{\sl Sufficiency} was proved in \cite{Omirov}.
\endproof
Let $E$ be an $n$-dimensional nilpotent evolution algebra with maximal nilpotent index. Then by the following scaling of basis

\begin{equation}\label{0}
\left\{\begin{aligned}
{}e'_1&=a_{12}^{-1/2} a_{23}^{-1/4}\cdots a_{n-1,n}^{-(1/2)^{n-1}}e_1\\
e'_2 &=a_{23}^{-1/2} a_{34}^{-1/4}\cdots a_{n-1,n}^{-(1/2)^{n-2}}e_2\\
& \dots\\
e'_{n-1}& =a_{n-1,n}^{-1/2}e_{n-1}\\
e'_n&=e_n
\end{aligned}\right.
\end{equation}
the evolution algebra is isomorphic to an evolution algebra $E'$
with matrix of structural constants
\[A'=\begin{pmatrix}
0 & 1 & a'_{13} &\dots &a'_{1n} \\
 0 & 0 & 1 &\dots &a'_{2n} \\
0 & 0 & 0 &\dots &a'_{3n} \\
\vdots & \vdots & \vdots &\cdots & \vdots \\
0 & 0 & 0 &\cdots & 1\\
0 & 0 & 0 &\cdots &0
\end{pmatrix}
 \]

Let $E$ be an $n$-dimensional nilpotent evolution algebra such that the matrix of structural constants  satisfies $a_{i_1i_1+1}=\dots=a_{i_si_s+1}=0$, for some $s=1,\dots,n-1$. Then omitting all multipliers $a_{i_ki_k+1}$, $k=1,\dots,s$
in \eqref{0} one can show that the evolution algebra $E$ is isomorphic to an evolution algebra
$E'$ with matrix of structural constants $A'=(a'_{ij})$, with $a'_{i_1i_1+1}=\dots=a'_{i_si_s+1}=0$ and
$a'_{ii+1}=1, \ \ i\ne i_1,\dots ,i_s$.

The following theorem gives the classification of evolution algebras with matrix of structural constants as $A'$.

\begin{thm}\label{th3}
Any finite-dimensional complex evolution algebra with maximal nilpotent index is isomorphic to one
of pairwise non-isomorphic algebras with matrix of structural constants
\[\begin{pmatrix}
 0 & 1 & a_{13} &\dots &a_{1,n-1}&0 \\[1.5mm]
 0 & 0 & 1 &\dots &a_{2,n-1}& 0 \\[1.5mm]
0 & 0 & 0 &\dots &a_{3,n-1}& 0 \\[1.5mm]
\vdots & \vdots & \vdots &\cdots & \vdots &\vdots\\[1.5mm]
0 & 0 & 0 &\cdots & 0& 1\\[1.5mm]
0 & 0 & 0 &\cdots &0 & 0
\end{pmatrix},
 \]
where one of non-zero $a_{ij}$ can be chosen equal to 1.
\end{thm}
\proof Assume $\varphi=(\alpha_{ij})_{i,j=1,\dots,n}$ is an isomorphism between evolution algebras $E$ and $E'$ with multiplication tables
\[E:\left\{
\begin{aligned}
{}e_i^2 &= e_{i+1}+\sum_{j=i+2}^na_{ij}e_j, && i=1,\dots, n-2,\\
e_{n-1}^2 &= e_n, &&\\
e_n^2 &= 0. &&
\end{aligned}
\right.\]
\[E':\left\{
\begin{aligned}
{}{e'_i}^2 &= e'_{i+1}+\sum_{j=i+2}^nb_{ij}e'_j, \ \ i=1,\dots, n-2,\\
{e'}^2_{\!\! n-1}  &= e'_n,\\
{e'}^2_{\! \! n} &= 0.
\end{aligned}
\right.\]

 We shall use the following lemma.

\begin{lemma}
$\varphi=(\alpha_{ij})$ is defined by
\begin{equation}\label{ll1}
\alpha_{ii}=\alpha_{11}^{2^{i-1}}, \ \ \alpha_{11}\ne 0; \qquad \alpha_{ij}=0, \ \ i\ne j,\, \ \ j\ne n; \quad  \alpha_{in}\in \mathbb C.
\end{equation}
\end{lemma}
\proof For $i=n$ consider
\[0=\varphi(e_n^2)={e'}^2_{\! \! n} =\sum_{j=1}^{n-2}\alpha_{nj}^2\Big(e_{j+1}+\sum_{s=j+2}^n a_{js}e_s\Big)+\alpha^2_{n,n-1}e_n.\]
From this equality it follows that $\alpha_{nj}=0$, $j=1,\dots, n-1$.

For $i=n-1$ we have
\[e_n'=\alpha_{nn}e_n=\varphi(e_{n-1}^2)=\sum_{i=1}^{n-1}\alpha_{n-1,i}^2e_i^2=\sum_{i=1}^{n-2}\alpha_{n-1,i}^2\Big(e_{i+1}+
\sum_{j=i+2}^n a_{ij}e_j\Big)+\alpha^2_{n-1,n-1}e_n.\]
From this equalities we get
 \[\alpha_{n-1,i}=0, \quad i=1,\dots, n-2, \qquad \alpha_{nn}=\alpha^2_{n-1, n-1}.\]
 Similarly from $e_{n-1}'e_j=0$, we get $\alpha_{j,n-1}=0$. Hence $\alpha_{n-1,j}=\alpha_{j,n-1}=0$, $j\ne n-1$.
 Assume \eqref{ll1} is true for $i\geq k$ and $j\ne n$. We shall check it for $i=k-1$.
 By the assumptions we have
  \[\varphi(e_{k-1}^2)=\sum_{s=1}^{k-1}\alpha_{k-1,s}^2e_s^2=\sum_{s=1}^{k-1}\alpha_{k-1,s}^2 \, \Big(e_{s+1}+
\sum_{j=s+2}^n a_{sj}e_j\Big).\]
 On the other hand

 \begin{align*}
\varphi(e_{k-1}^2)&={e'}^2_{\! \! k-1} =e_k'+\sum_{j=k+1}^nb_{k-1,j}e'_j \\
{} & = \alpha_{kk}e_k+\alpha_{kn}e_n+\sum_{j=k+1}^{n}b_{k-1,j} \, (\alpha_{jj} \, e_j+\alpha_{jn} \, e_n).
\end{align*}
Hence $\alpha_{k-1,j}=0$, $j=1,\dots,k-2$. Using the above equalities we get
\begin{equation}\label{lk}
\alpha_{k-1,k-1}^2e_k+\Big(\sum_{j=k+1}^na_{k-1,j}e_j\Big)\alpha_{k-1,k-1}^2=
\alpha_{kk}e_k+\alpha_{kn}e_n+\sum_{j=k+1}^{n}b_{k-1,j}\,(\alpha_{jj}e_j+\alpha_{jn}e_n)
\end{equation}
which gives $\alpha_{kk}=\alpha^2_{k-1, k-1}$.

For $s<k$ using assumptions,  we get
\begin{align*}
0&=\varphi(e_{k-1}e_s)=e_{k-1}'e_s'=\big(\alpha_{k-1,k-1}e_{k-1}+\alpha_{k-1,n}e_n\big)
\Big(\sum_{t=1}^k\alpha_{st}e_t+\alpha_{sn}e_n\Big)\\
{}& =\alpha_{k-1,k-1}\alpha_{s,k-1}e^2_{k-1}=\alpha_{k-1,k-1}\alpha_{s,k-1}\Big(e_k+\sum_{j=k+1}^na_{kj}e_j\Big).
\end{align*}
Then $\alpha_{s,k-1}=0$. Hence
\[\alpha_{j,k-1}=0, \ \  j<k-1, \qquad \alpha_{k-1,j}=0, \ \ j\ne k-1. \qedhere \]
\endproof
Now we shall continue the proof of theorem. From the equality \eqref{lk} we get
\begin{multline*}
\alpha_{k-1,k-1}^2  \, a_{k-1,n}  \, e_n+\Big(\sum_{j=k+1}^{n-1}a_{k-1,j}  \, e_j\Big)\alpha_{k-1,k-1}^2 \\
=\sum_{j=k+1}^{n-1}b_{k-1,j}  \, \alpha_{jj} \, e_j+\Big(\alpha_{kn}+2b_{k-1,n} \, \alpha_{nn}+\sum_{j=k+1}^{n-1}b_{k-1,j}\alpha_{jn}\Big)e_n.
\end{multline*}

Consequently,
\begin{align*}
\alpha_{k-1,k-1}^2a_{k-1,n}&=\alpha_{kn}+2b_{k-1,n} \, \alpha_{nn} +\sum_{j=k+1}^{n-1}b_{k-1,j}\alpha_{jn}, && k=2,\dots,n-1. \\
\alpha_{k-1,k-1}^2a_{k-1,j}&=\alpha_{jj}b_{k-1,j}, && j=k+1,\dots,n-1.
\end{align*}
From these formulas using $\alpha_{kk}=\alpha_{11}^{2^{k-1}}$, we obtain
\begin{equation}
\begin{aligned}\label{ss}
\alpha_{11}^{2^{k-1}}a_{k-1,n}&=\alpha_{kn}+2b_{k-1,n} \, \alpha_{11}^{2^{n-1}}+\sum_{j=k+1}^{n-1}b_{k-1,j}\alpha_{jn}, && k=2,\dots,n-1. \\
\alpha_{11}^{2^{k-1}}a_{k-1,j}&=\alpha_{11}^{2^{j-1}}b_{k-1,j},&& j=k+1,\dots,n-1.
\end{aligned}
\end{equation}

From the second equation of the system \eqref{ss} we obtain
\begin{equation}\label{bu}
b_{k-1,j}=\alpha_{11}^{2^{k-1}-2^{j-1}}a_{k-1,j}, \qquad k=2,\dots, n-1, \ \ j=k+1,\dots,n-1.
\end{equation}

Using (\ref{bu}), in order to have $b_{k-1,n}=0$, in the first equation of \eqref{ss} we put
$$\alpha_{kn}=\alpha_{11}^{2^{k-1}}a_{k-1,n}- \sum_{j=k+1}^{n-1}b_{k-1,j}\alpha_{jn}$$ $$
=\alpha_{11}^{2^{k-1}}a_{k-1,n}- \sum_{j=k+1}^{n-1}\alpha_{11}^{2^{k-1}-2^{j-1}}a_{k-1,j}\alpha_{jn},\ \ k=2,\dots,n-1.$$

If there exist $k_0$, $j_0$ such that $a_{k_0,j_0}\ne 0$ then taking $\alpha_{11}= a_{k_0,j_0}^{\frac{-1} {2^{k_0}-2^{j_0-1}}}$, we have
\[b_{k_0,j_0}=1, \quad k_0=1,\dots,n-2. \qedhere \]
\endproof

\begin{rk} Note that the evolution algebras of Theorem \ref{th3} are also algebras
of maximal index of solvability and maximal index of right nilpotency.
\end{rk}
\begin{ex} Any four-dimensional complex evolution algebra with maximal nilpotent index is
isomorphic to one of the algebras with the following matrix of structural constants
\[\begin{pmatrix}
0&1&0&0\\
0&0&1&0\\
0&0&0&1\\
0&0&0&0
\end{pmatrix},  \quad
\begin{pmatrix}
0&1&1&0\\
0&0&1&0\\
0&0&0&1\\
0&0&0&0
\end{pmatrix}.\]
Any five-dimensional complex evolution algebra with  maximal nilpotent index is
isomorphic to one of the algebras with the following matrix of structural constants
\[\begin{pmatrix}
0&1&0&0&0\\
0&0&1&0&0\\
0&0&0&1&0\\
0&0&0&0&1\\
0&0&0&0&0
\end{pmatrix},  \ \ \begin{pmatrix}
0&1&0&0&0\\
0&0&1&1&0\\
0&0&0&1&0\\
0&0&0&0&1\\
0&0&0&0&0
\end{pmatrix},  \ \ \begin{pmatrix}
0&1&0&1&0\\
0&0&1&d&0\\
0&0&0&1&0\\
0&0&0&0&1\\
0&0&0&0&0
\end{pmatrix}, \ \ \begin{pmatrix}
0&1&1&b&0\\
0&0&1&d&0\\
0&0&0&1&0\\
0&0&0&0&1\\
0&0&0&0&0
\end{pmatrix},\] where $b,d\in \mathbb{C}$.
\end{ex}
\begin{proposition} Let $E$ be an evolution algebra with matrix of structural constants
\[A=\begin{pmatrix}
 0 & 1 & a_{13} &\dots & a_{1,m+1}& a_{1,m+2} &\dots &a_{1,n-1}&a_{1n} \\[1.5mm]
 0 & 0 & 1 &\dots & a_{2,m+1} & a_{2,m+2} &\dots &a_{2,n-1}& a_{2n} \\[1.5mm]
0 & 0 & 0 &\dots & a_{3,m+1} & a_{3,m+2} &\dots &a_{3,n-1}& a_{3n} \\[1.5mm]
\vdots & \vdots & \vdots &\vdots &\vdots &  \vdots &\vdots & \vdots &\vdots\\[1.5mm]
0 & 0 & 0 &\dots & a_{m-1,m+1}  & a_{m-1,m+2} &\dots &  a_{m-1,n-1} & a_{m-1,n}\\[1.5mm]
0 & 0 & 0 &\dots & a_{m,m+1}  & a_{m,m+2} &\dots &  a_{m,n-1} & a_{mn}\\[1.5mm]
\vdots & \vdots & \vdots &\vdots &\vdots & \vdots &\vdots & \vdots &\vdots\\[1.5mm]
0 & 0 & 0 &\dots &0  &0& \dots & 0& 1\\[1.5mm]
0 & 0 & 0 &\dots &0  &0 & \dots &0 & 0
\end{pmatrix},\]
with $a_{i,i+1}=1$, for $i\ne m$,  \ $i=1,\dots, n-1$, \
$a_{m,m+1}=0$ and $a_{m,m+2}  + a_{m-1,m+1}^2   \neq 0.$
 Then
\begin{equation}\label{k1}
E^{2^k+1}=E^{2^k+2}=\dots=E^{2^{k+1}}=\left\{\begin{array}{l}
\langle e^2_{k+1}, \dots, e^2_{m-1}, e_{m+2}, \dots, e_n\rangle, \quad \text{if} \ \ k\leq m-2,\\[2mm]
\langle e_{k+3}, \dots, e_n\rangle, \  \qquad \qquad \text{if} \ \ m-1\leq k\leq n-3,
\end{array}
\right.
\end{equation}
and $E^{2^{n-2}+1}=0$, i.e., its nilpotent index is $2^{n-2}+1$.
\end{proposition}
\proof  We have
\[e^2_i=e_{i+1}+\sum_{j=i+2}^n a_{ij}e_j, \ \ i=1,\dots,m-1, m+1,\dots,n-1; \qquad e^2_m=\sum_{j=m+2}^n a_{mj}e_j. \]
\[E^2=\langle e_1^2, e^2_2, \dots, e^2_{n-1}\rangle.\]
It is easy to see that $\langle e_i^2, \, i=m, m+1,\dots,n\rangle
= \langle e_i, \, i=m+2,\dots,n\rangle$ and $e_1^2, e^2_2, \dots,
e^2_{m-1}$, $e_{m+2},\dots, e_n$ are linearly independent. Thus
\[E^2=\langle e_1^2, e^2_2, \dots, e^2_{m-1}, e_{m+2},\dots, e_n\rangle.\]
\[E^3=EE^2=\langle e^2_ie_j, \, e^2_k \mid i=1,\dots,m-1, \ j=i+1,\dots,n-1, \ k=m+2, \dots, n-1\rangle.\]
We have
\[
e_i^2e_j=
\begin{cases}
e_j^2, &\text{if \ $j=i+1$,}\\
a_{ij} \, e_j^2, &\text{if \ $j>i+1$.}
\end{cases}
\]

In case $a_{m,m+2}\ne 0$, from $\displaystyle e_m^2=a_{m,m+2} \, e_{m+2}+\sum_{j=m+3}^na_{mj} \, e_j$,
we obtain $e_{m+2}\in E^3$. If $a_{m-1,m+1}\ne 0$, then since $e^2_{m-1} \, e_{m+1}=a_{m-1,m+1} \, e^2_{m+1}$, we
conclude that $e_{m+2}\in E^3$. Thus we obtain
\[E^3=\langle e_2^2, \dots, e^2_{m-1}, e_{m+2},\dots, e_n\rangle.\]
Now we shall compute $E^4=EE^3+E^2E^2$. Similarly as in the case $EE^2$, we get
\begin{align*}
EE^3&=\langle e^2_ie_j,\, e^2_k \mid i=2,\dots,m-1, \ j=i+1,\dots,n-1, \
k=m+2, \dots, n-1\rangle \\
{}&=\langle e_3^2, \dots, e^2_{m-1}, e_{m+2},\dots, e_n\rangle. \\
E^2E^2 &=\langle e^2_i, \, e^2_je_k, \, e^2_pe^2_q \mid i,k=m+2,\dots, n-1, \ j, p,q=1,\dots,m-1\rangle.
\end{align*}
It is easy to see that
\[e^2_ie_k=a_{ik}e_k^2, \qquad i=1,\dots,m-1, \ \ k=m+2,\dots,n-1.\]
\[
e^2_pe^2_q=
\begin{cases}
e^2_{p+1}, &\text{if \ $p=q=1,\dots,m-1$,}\\
\displaystyle a_{p,q+1}e^2_{q+1}+\sum_{j=q+2}^na_{pj}a_{qj}e_j^2, &\text{if \ $p<q$.}
\end{cases}
\]
Using these equalities, we obtain
\[E^2E^2=\langle e_2^2, \dots, e^2_{m-1}, e_{m+2},\dots, e_n\rangle=E^3.\]
Consequently, $E^4=E^3$.

Lets assume that the equalities \eqref{k1} are true for $k$, we shall prove it for $k+1$.
\begin{align*}
E^{2^{k+1}+1}&=EE^{2^{k+1}}+E^2E^{2^{k+1}-1}+\dots+E^{2^k}E^{2^k+1}=\big(E+E^2+\dots+E^{2^k}\big)E^{2^k+1}\\
{}&=EE^{2^k+1}=\langle e_1, \dots, e_n\rangle\langle e^2_{k+1}, \dots, e^2_{m-1}, e_{m+2}, \dots, e_n\rangle \\
{}&=\langle e^2_{k+2}, \dots, e^2_{m-1}, e_{m+2}, \dots, e_n\rangle.
\end{align*}
We also have
\begin{align*}
E^{2^{k+2}}&=EE^{2^{k+2}-1}+E^2E^{2^{k+2}-2}+\dots+E^{2^{k+1}}E^{2^{k+1}}\supseteq E^{2^{k+1}}E^{2^{k+1}}\\
{}&=\langle e^2_{k+2}, \dots, e^2_{m-1}, e_{m+2}, \dots, e_n\rangle.
\end{align*}
Moreover, we obtain
\begin{align*}
\langle e^2_{k+2}, \dots, e^2_{m-1}, e_{m+2}, \dots, e_n\rangle& =E^{2^{k+1}+1}\supseteq E^{2^{k+1}+2}\supseteq\dots\supseteq E^{2^{k+2}}\\
{}& \supseteq \langle e^2_{k+2}, \dots, e^2_{m-1}, e_{m+2}, \dots, e_n\rangle.
\end{align*}
Consequently
\[E^{2^{k+1}+1}= E^{2^{k+1}+2}=\dots =
E^{2^{k+2}}=\langle e^2_{k+2}, \dots, e^2_{m-1}, e_{m+2}, \dots,
e_n\rangle.\] This gives the formula \eqref{k1} for $k \leq m-2.$

\

Now we will show for $k=m-1.$  Since the formula \eqref{k1} is
true for $k = m-2,$ we have
$$E^{2^{m-2}+1}=E^{2^{m-2}+2}=\dots = E^{2^{m-1}}=\langle e_{m-1}^2, e_{m+2}, \dots, e_n, \rangle.$$

Consider $$E^{2^{m-1}+1} = E E^{2^{m-1}} + E^2 E^{2^{m-1}-1} +
\dots +E^{2^{m-2}} E^{2^{m-2}+1}.$$

Since $$E E^{2^{m-1}} = \langle a_{m-1,m+1}e_{m+2},
a_{m,m+2}e_{m+2}, e_{m+3}, \dots,e_n \rangle$$ and $a_{m,m+2}  +
a_{m-1,m+1}^2   \neq 0$ we have
$$E^{2^{m-1}+1} = \langle e_{m+2}, e_{m+3}, \dots,e_n \rangle.$$

Now consider  $$E^{2^{m-1}+2} = E E^{2^{m-1}+1} + E^2 E^{2^{m-1}}
+ \dots +E^{2^{m-2}+1} E^{2^{m-2}+1}.$$

Taking into account that $$E^{2^{m-2}+1} E^{2^{m-2}+1} = \langle e_{m-1}^2e_{m-1}^2,
e_{m+3}, \dots,e_n \rangle = \langle (a_{m,m+2}  +
a_{m-1,m+1}^2)e_{m+2}, e_{m+3}, \dots,e_n \rangle $$ and 
$a_{m,m+2} + a_{m-1,m+1}^2 \neq 0,$ we conclude 
$$E^{2^{m-1}+2} = \langle e_{m+2}, e_{m+3}, \dots,e_n \rangle.$$

Analogously, we obtain
$$E^{2^{m-1}+3} = E^{2^{m-1}+4} = \dots = E^{2^{m}} = \langle e_{m+2}, e_{m+3}, \dots,e_n
\rangle.$$ Thus we obtain the restriction for $k=m-1.$

\

For $k \geq m$ is easy to get the formula \eqref{k1} and using the
formula, we get $E^{2^{n-2}+1}=0$.
\endproof

The following example shows that the nilpotent index may be $2^{n-s}+1$ for any $s=1,\dots,n-1$.

\begin{ex}\label{ex2} Consider an $n$-dimensional nilpotent evolution algebra
$E$ with  matrix of structural constants $A=(a_{ij})$ satisfying
\begin{align*}
a_{i_kj}&=0, && j=1,2,\dots,n-1, && k=1,\dots,s-1;\\
 a_{ji_k}&=0, && j=1,2,\dots,n, && k=1,\dots,s-1,
\end{align*}
where $1\leq i_1< i_2<\dots< i_{s-1}\leq n-2$, $s<n$.
Then the nilpotent index of $E$ is $2^{n-s}+1$.
\end{ex}
The following proposition generalizes Example \ref{ex2}.

\begin{proposition} Let $E$ be an $n$-dimensional evolution algebra with matrix of structural constants
$A=(a_{ij})$ such that for some $s<n-1$ and $1\leq i_1<i_2<\dots<i_{s-1}\leq n-2$ satisfies

\begin{align*}
a_{ji_k}=a_{i_kj}&=0, && \text{for all} \ \ j\notin\{i_1,\dots,i_{s-1},n\}, && k=1,\dots, s-1, \\
a_{i_ki_{k+1}}&=1,&& \text{for all} \ \ k=1,\dots, s-1. &&
\end{align*}
i.e. $A$ has the following form
{\small
\[\setcounter{MaxMatrixCols}{13}
\begin{pmatrix}
0&1&a_{13}&\dots&a_{1i_1-1}&0&a_{1i_1+1}&\dots&a_{1i_{s-1}-1}&0&a_{1i_{s-1}+1}&\dots&a_{1n}\\[2mm]
0&0&1&\dots&a_{2i_1-1}&0&a_{2i_1+1}&\dots&a_{2i_{s-1}-1}&0&a_{2i_{s-1}+1}&\dots&a_{2n}\\[2mm]
\vdots&\vdots&\vdots&\dots&\vdots&\vdots&\vdots&\dots&\vdots&\vdots&\vdots&\dots&\vdots\\[2mm]
0&0&0&\dots&0&1&a_{i_1-1i_1+1}&\dots&a_{i_1-1i_{s-1}-1}&0&a_{i_1-1i_{s-1}+1}&\dots&a_{i_1-1n}\\[2mm]
0&0&0&\dots&0&0&0&\dots&0&a_{i_1i_{s-1}}&0&\dots&a_{i_1n}\\[2mm]
0&0&0&\dots&0&0&0&\dots&a_{i_1+1i_{s-1}-1}&0&a_{i_1+1i_{s-1}+1}&\dots&a_{i_1+1n}\\[2mm]
\vdots&\vdots&\vdots&\dots&\vdots&\vdots&\vdots&\dots&\vdots&\vdots&\vdots&\dots&\vdots\\[2mm]
0&0&0&\dots&0&0&0&\dots&0&0&0&\dots&1\\[2mm]
0&0&0&\dots&0&0&0&\dots&0&0&0&\dots&0\\[2mm]
\end{pmatrix}.\]}
Then the nilpotent index of $E$ is equal to $2^{\max\{s-1, n-s\}}+1$.
\end{proposition}
\proof The evolution algebra $E$ can be written as $E=A+B$, where $A=\langle e_i \mid i\ne i_1,\dots, i_{s-1}\rangle$,
$B=\langle e_{i_1},e_{i_2},\dots, e_{i_{s-1}}\rangle$. It is easy to see that $AB=0$, this implies $A^iB^j=0$, $i,j=1,2,\dots$
Consequently, $E^k=A^k+B^k$. Using similar arguments as above (for computation of the maximal nilpotent index) one can see that
the nilpotent index of $A$ is $2^{n-s}+1$ and the nilpotent index of $B$ is $2^{s-1}+1$. This completes the proof.
\endproof

\begin{rk} By   \cite[Proposition 4.7]{Omirov} if $E$ is an $n$-dimensional nilpotent evolution algebra
with index of nilpotency not equal to $2^{n-1}+1$, then it is not greater than $2^{n-2}+1$.
Moreover, in the paper \cite{Omirov}, there is an example of evolution algebra
with nilpotent index $3\cdot 2^{k-4}+1$, where $4\leq k\leq n$.
Therefore it is interesting to know all possible values of the nilpotent index for nilpotent evolution algebras.
This problem is difficult, but for small values of $n$ one can do exact calculations.
For example, if $n=3$ then the nilpotent index can be 2,3,5.
 For $n=4$ all possible values of the nilpotent index are 2,3,4,5,9.
 The 4-dimensional evolution algebra $E$, with the following matrix
\[A=\begin{pmatrix}
 0&1&b&c\\
 0&0&0&-b^2f\\
 0&0&0&f\\
 0&0&0&0
 \end{pmatrix},\]
 has nilpotent index 4, where $bf\ne 0$. This case is interesting since it has not the form $2^k+1$.
\end{rk}
\section{Dibaric algebras} \label{S:dib}

In this section we will study  some dibaricity properties of arbitrary algebras and evolution algebras.

\begin{thm}\label{thm2} Any finite-dimensional nilpotent evolution algebra $E$ is not dibaric.
\end{thm}
\proof
Assume $\varphi=(b_{ij})_{i=1,\dots,n; \, j=1,2}$ is a homomorphism $\varphi \colon E \to \mathfrak A$.
We shall  use the following

\begin{lemma}\label{li} For any $i,j=1,\dots,n$, we have
\begin{equation}\label{el}
b_{i1}b_{j2}=b_{i2}b_{j1}=0.
\end{equation}
\end{lemma}
\proof Without loss of generality we assume $i\leq j$ and use mathematical induction. Let $C_k$ denote all cases
of \eqref{el} with $2n -(i+j)+1 = k$. For $k=1$, i.e. $i=j=n$, from $\varphi(e_n^2)=0$ we get
\begin{equation}\label{i1}
b_{n1}b_{n2}=0.
\end{equation}
For $k=2$ we have $i=n-1$ and $j=n$. We get
\[0=\varphi(e_{n-1}e_n)=\frac{1}{2}(b_{n-1,1}b_{n2}+b_{n-1,2}b_{n1})(m+w).\]
 This by \eqref{i1} gives
\[b_{n-1,1}b_{n2}=b_{n-1,2}b_{n1}=0.\]
Assuming that $C_k$ holds, we have to prove $C_{k+1}$, that is, equation \eqref{el}, for any
$i,j=1,\dots,n$, $i\leq j$, which satisfy $2n-(i+j)+1=k+1$.

{\it Case} $i<j$: From $\varphi(e_ie_j)=\frac{1}{2}(b_{i1}b_{j2}+b_{i2}b_{j1})(m+w)=0$ we get
\begin{equation}\label{i2}
b_{i1}b_{j2}+b_{i2}b_{j1}=0.
\end{equation}
By assumptions we have $i<j$, $2n-2j+1\leq k$ and $b_{j1}b_{j2}=0$. This by \eqref{i2} gives \eqref{el}.

{\it Case} $i=j$: Consider
\begin{align*}
\varphi(e_i^2)&=\varphi\Big(\sum_{s=i+1}^na_{is} \, e_s\Big)=\sum_{s=i+1}^n a_{is} \, \varphi(e_s)
=\Big(\sum_{s=i+1}^n a_{is}b_{s1}\Big)m+\Big(\sum_{s=i+1}^n a_{is}b_{s2}\Big)w\\{}&=\varphi(e_i)^2=\frac{1}{2} b_{i1}b_{i2}(m+w).
\end{align*}
Consequently,
\begin{equation}\label{i3}
\left\{\begin{aligned}
2\sum_{s=i+1}^n a_{is}b_{s1}&=b_{i1}b_{i2},\\
2\sum_{s=i+1}^n a_{is}b_{s2}&=b_{i1}b_{i2}.
\end{aligned}\right.
\end{equation}
Since $2n-(i+s)+1\leq k$ for any $s=i+1, i+2, \dots, n$, by the assumption of the induction we get
\begin{equation}\label{i4}
b_{s1}b_{i2}=b_{s2}b_{i1}=0.
\end{equation}
Now, multiplying  both sides of the first equation of \eqref{i3} by $b_{i2}$, then by \eqref{i4} we get
$b_{i1}b_{i2}=0$.
\endproof

 Now we shall continue the proof of the theorem.  By Lemma \ref{li}, if there exists $i_0$ such that $b_{i_01}\ne 0$ then $b_{j2}=0$ for all $j$, i.e., $\varphi(e_i)=b_{i1}m$. Such $\varphi$ is not onto.
\endproof

The following result gives a sufficient condition for an arbitrary algebra to be non dibaric.

\begin{thm}  Let $\mathbf A$ be a finite-dimensional real algebra with table of multiplication $\displaystyle e_ie_j=\sum_k a^k_{ij}e_k$, where $(a^k_{ij})_{i,j,k=1,\dots,n}$  is the matrix of structural constants, and such that
the matrix $A=(a^k_{ii})_{i,k=1,\dots,n}$ has $\det(A)\ne 0$. Then $\mathbf A$ is not dibaric.
\end{thm}
\proof
Assume $\varphi=(\alpha_{ij})_{i=1,\dots,n; \, j=1,2}$ is a homomorphism $\varphi \colon \mathbf A \to \mathfrak A$.
We have
\begin{align*}
\varphi(e^2_i)&=\sum^n_{s=1}a_{ii}^s  \, (\alpha_{s1}m+\alpha_{s2}w). \\
\varphi(e^2_i)&=\frac{1}{2}(\alpha_{i1}\alpha_{i2})(m+w).
\end{align*}
Consequently,
\begin{equation}\label{u1}
\left\{\begin{aligned}
2\sum_{s=1}^n a_{ii}^s\alpha_{s1}&=\alpha_{i1}\alpha_{i2},\\
2\sum_{s=1}^n a_{ii}^s\alpha_{s2}&=\alpha_{i1}\alpha_{i2}.
\end{aligned}\right.
\end{equation}
Subtracting from first equation of the system \eqref{u1} the second one, we obtain
\begin{equation}\label{u2}
\sum_{s=1}^n a_{ii}^s(\alpha_{s1}-\alpha_{s2})=0, \quad i=1,\dots,n.
\end{equation}
 If $\det(A)\ne 0$ we get from the system \eqref{u2} that $\alpha_{i1}=\alpha_{i2}$ for all $i$.
 Hence $\varphi(e_i)=\alpha_{i1}(m+w)$, but such $\varphi$ is not onto.
 \endproof
 \begin{rk}Non dibaric algebras given by Theorem \ref{thm2} show that the
 condition $\det(A)\ne 0$ is not necessary to be non dibaric.
\end{rk}
 \begin{cor} Let $E$ be an evolution algebra with matrix $A$ of structural constants. If $\det(A)\ne 0$ then $E$ is not dibaric.
\end{cor}

\begin{defn}[\cite{LOR}] For a given algebra $\mathbf A$, a pair $(f,g)$ of linear forms $f \colon \mathbf A\to \R$,
$g \colon \mathbf A\to \R$ is called bq-homomorphism if
\begin{equation}\label{bq}
f(xy)=g(xy)=\frac{f(x)g(y)+f(y)g(x)}{2}\ \ \text{for any} \ \ x, y\in \mathbf A.
\end{equation}
\end{defn}
Note that if $f=g$ then  condition \eqref{bq} implies that $f$ is a homomorphism.

A  bq-homomorphism $(f,g)$ is called non-zero if both $f$ and $g$ are non-zero.

\begin{thm}[\cite{LOR}]\label{t0}  An algebra $\mathbf A$ is dibaric if and only if there is a non-zero bq-homomorphism $(f,g)$.
\end{thm}
In case of Theorem \ref{t0} the homomorphism
$\varphi \colon \mathbf A\to \mathfrak A$ has the form $\varphi(x)=f(x)m+g(x)w$. Lets denote
\[V_n=\{x\in \mathbf A: \varphi(x^n)=0\}.\]
\begin{proposition} For any $n\geq 3$, we have $V_n=V_3$.
\end{proposition}
\proof  Using $(m+w)^n=m+w$ and mathematical induction, one can prove the following formula
\[\varphi(x^n)=\frac{f(x)g(x)}{2^{n-1}}\big(f(x)+g(x)\big)^{n-2}(m+w).\]
From this formula, for $n\geq 3$, we get $\varphi(x^n)=0$ if and only if $f(x)g(x)\big(f(x)+g(x)\big)=0$. This completes the proof.
\endproof
\begin{rk} Any solvable algebra $\mathbf A$ is not dibaric. Indeed,  there is an homomorphism $\varphi$ onto $\mathfrak A$.
It is easy to see that
\begin{equation}\label{sb}
\varphi(\mathbf A^{[k]})=\varphi(\mathbf A)^{[k]}=
\mathfrak A^{[k]}=\mathfrak A^2\cong \mathbb R, \quad \text{ for all } \, k\geq 2.
 \end{equation} By the solvability of $\mathbf A$ there exists $k$ such that $\mathbf A^{[k]}=0$.
 Then from \eqref{sb} we get $0=\mathfrak A^2\cong \mathbb R$, this is a contradiction.
\end{rk}

 \subsection*{Two-dimensional dibaric evolution algebras}
 In this subsection we find a criterion for two-dimensional real  evolution algebra to be dibaric.

 Let the two-dimensional real evolution algebra $E$ be given by the matrix of structural constants
 $A=\begin{pmatrix}
 a&b\\
 c&d
 \end{pmatrix}$.

\begin{proposition}
The two-dimensional real evolution algebra $E$ is dibaric if and only if one of the following conditions hold
\begin{enumerate}
  \item $b=d=0$ and $ac<0$;
  \item $b\ne 0$, $ad=bc$, $D\geq 0$ and $B^2+C^2\ne 0$,
\end{enumerate}
where $D=(8a-1)^2-32(bd+a^2)$,  $B=4a^2+4bd-a+a\sqrt{D}$ and $C=4a^2+4bd-a-a\sqrt{D}$.
\end{proposition}
\proof Let $\varphi=\begin{pmatrix}
 \alpha & \beta\\
 \gamma & \delta
  \end{pmatrix}$ be a homomorphism. It is onto if $\alpha\delta\ne \gamma\beta$. Moreover, $\varphi$ must satisfy
 \begin{equation}\label{s}
 \left\{\begin{aligned}
 2(a\alpha+b\gamma)&=\alpha\beta\\
 2(a\beta+b\delta)&=\alpha\beta\\
 2(c\alpha+d\gamma)&=\gamma\delta\\
 2(c\beta+d\delta)&=\gamma\delta\\
 \alpha\delta+\beta\gamma &=0
 \end{aligned}
 \right.
 \end{equation}
 The proof follows from the elementary analysis of the system \eqref{s}.
 \endproof

\section*{ Acknowledgements}

 The first and second authors were supported by Ministerio
de Ciencia e Innovaci\'on (European FEDER support included), grant
MTM2009-14464-C02, and by Xunta de Galicia, grant Incite09 207 215
PR.  The third and fourth authors
thank the Department of Algebra, University of Santiago de
Compostela, Spain,  for providing financial support to their visit to
the Department.
We thank J. Carlos Guti\'errez Fern\'andez for the carefully reviewing of our paper \cite{LAA}.


\begin{thebibliography}{20}

\bibitem{Omirov} Camacho L.M., G\'omez J.R., Omirov B.A.,
Turdibaev R.M. \emph{Some properties of evolution algebras}, Bull. Korean Math. Soc. 50(5) (2013) 1481--1494.

\bibitem{CG} Camacho L.M., G\'omez J.R., Omirov B.A., Turdibaev R.M. \emph{The derivations of some evolution algebras}, Linear Multilinear Algebra 61(3) (2013) 309--322.

\bibitem{LAA}  Casas J.M.,  Ladra M., Omirov B.A.,  Rozikov U.A., \emph{On nilpotent index and dibaricity of evolution algebras}.  Linear Algebra Appl. 439(1) (2013) 90--105.

\bibitem{CLOR} Casas J.M., Ladra M., Omirov B.A.,  Rozikov U.A.  \emph{On evolution algebras},   Algebra Colloquium. 21(2) (2014) 331--342.

\bibitem{LOR} Ladra M., Omirov B.A.,  Rozikov U.A. \emph{ Dibaric and evolution algebras in biology}.  Lobachevskii  J. Math. 35(3) (2014) 198--210.

\bibitem{ly} Lyubich Y.I. \emph{Mathematical structures in population
genetics}, Springer-Verlag, Berlin, 1992.

\bibitem{M} Mazzola G. \emph{The algebraic and geometric classification
of associative algebras of dimension five}, Manuscripta Math. 27
(1979) 81--101.

\bibitem{m} Reed M.L. \emph{Algebraic structure of genetic inheritance},  Bull. Amer. Math.
Soc. (N.S.)  34(2)  (1997) 107--130.

\bibitem{Tian2} Rozikov U.A., Tian J.P. \emph{Evolution algebras
generated by Gibbs measures}, Lobachevskii J. Math. 32(4) (2011) 270--277.

\bibitem{Tian1} Tian J.P. \emph{Evolution algebras and their
applications}, Lecture Notes in Math. 1921, Springer-Verlag, Berlin, 2008.

\bibitem{Umlauf} Umlauf K.A. \emph{\"Uber die Zusammensetzung der
endlichen continuierlichen transformationsgruppen insbesondere der
Gruppen vom Range null}, Thesis, Universit\"{a}t Leipzig, 1891.

\bibitem{w} W\"orz-Busekros A. \emph{Algebras in genetics}, Lecture Notes in
Biomathematics 36, Springer-Verlag, Berlin-New York, 1980.

\end{thebibliography}
\end{document}